\documentclass[11pt,reqno]{amsart}

\usepackage{enumerate,url,amssymb,  mathrsfs}%, pdfsync}% ,showkeys, pdfsync}
\usepackage{graphicx}
\newtheorem{theorem}{Theorem}[section]

\newtheorem*{lemma*}{Lemma}

\theoremstyle{definition}

\theoremstyle{remark}
\newtheorem{remark}[theorem]{Remark}

\numberwithin{equation}{section}

\newcommand{\D}{\mathbb{D}}

\def\XXint#1#2#3{{\setbox0=\hbox{$#1{#2#3}{\int}$}
\vcenter{\hbox{$#2#3$}}\kern-.5\wd0}}

\def\ge{\geqslant}
\newcommand{\E}{\mathcal{E}}

\begin{document}
\baselineskip6mm
\vskip0.4cm
\title[On a H\"older constant in quasiconformal mappings]{On a H\"older constant in the theory of quasiconformal mappings}

\
\author[I. Prause]{Istv\'an Prause}
\address{University of Helsinki, Department of Mathematics and Statistics,
         P.O. Box 68 , FIN-00014 University of Helsinki, Finland}
\email{istvan.prause@helsinki.fi}

\dedicatory{Dedicated to the memory of F.~W.~Gehring}

%\footnotetext[1]{submitted to F.~W.~Gehring Memorial Volume,
%Computational Methods and Function Theory}

\thanks{ The author was supported by the Academy of Finland grant 1266182.}


%\date{\today}

%\keywords{Sharp Bounds and Extremal Quasiconformal Mappings}
\maketitle

\begin{abstract} 
A $K$-quasiconformal selfmap of the unit disk with identity boundary values satisfies the H\"older estimate
$|f(z)-f(w)| \leq 4^{1-1/K} |z-w|^{1/K}$. The constant $4^{1-\frac{1}{K}}$ is sharp.
\end{abstract}

\section{Introduction}
A classical result in the theory of quasiconformal mappings, known as Mori's theorem is the following. If $f \colon \D \to \D=f(\D)$, $f(0)=0$ is $K$-quasiconformal map of the unit disk $\mathbb{D}=B(0,1)$ then
\begin{equation}
\label{eq:mori}
 |f(z)-f(w)| \leq 16 |z-w|^{1/K} \quad z,w \in \mathbb{D}.
\end{equation}
See for instance \cite{mori,ahlfors,lehto-virtanen}. Here the constant $16$ is optimal as an absolute constant, however it has been conjectured in \cite[p.~68]{lehto-virtanen} that $16$ could be replaced by $16^{1-1/K}$ if we allow dependence on $K$. We refer to the texts \cite{ahlfors,lehto-virtanen,aim09} for different definitions and basic properties of quasiconformal mappings.

The purpose of this note is to point out the following sharp counterpart of \eqref{eq:mori}. Below, we require identity boundary values, in which case the requirement $f(0)=0$ may be omitted.

\begin{theorem}
Let $f \colon \D \to \D$ be a $K$-quasiconformal mapping with (boundary extension) $f(z)=z$ for $|z|=1$. Then for every $z,w \in \D$ we have
\begin{equation} 
\label{eq:main}
|f(z)-f(w)| \leq 4^{1-1/K} |z-w|^{1/K}.
\end{equation}
Moreover, $4^{1-1/K}$ cannot be replaced by any smaller number depending only on $K$.
\end{theorem}

A theorem of this type in $\mathbb{R}^n$, $n \geqslant 2$ has recently appeared in \cite{vuorinen-zhang}. Their constant for $n=2$ gives the non-optimal value
\[ 4^{1-1/K} \cdot 2^{1-1/K} K^{\frac{1}{2K}} \left( \frac{K}{K-1} \right)^{\frac{1-1/K}{2}}.
\]

\begin{remark}
As it will be clear from the proof, the same bound holds for any $K$-quasiconformal principal mapping conformal outside the unit disk.
\end{remark}

\section{Proof}

\begin{proof}
Let us extend $f$ identically outside the unit disk. This way $f$ becomes a global quasiconformal map of the complex plane with Beltrami coefficient $\mu$, where $\|\mu\|_\infty \leq \frac{K-1}{K+1}$.
We embed $f$ in a holomorphic flow of solutions $\{f^\lambda\}_{\lambda \in \D}$ as follows.
For $\lambda \in \mathbb{D}$ we solve the uniformly elliptic Beltrami equation
\begin{equation}
\label{eq:beltrami} 
f_{\bar z}^{\lambda}= \lambda \mu \cdot \frac{K+1}{K-1} f_z^{\lambda},
\end{equation}
under the normalization of the so-called principal solution. As $\mu$ vanishes outside the unit disk, $f^\lambda$ will be conformal outside $\D$ and we require the asymptotics $f(z)=z +o(1)$ as $z \to \infty$.
According to a classical Koebe-type distortion theorem \cite[Theorem 1.4]{pommerenke}
\begin{equation} 
\label{eq:koebe}
|f^\lambda(z)| <2, \quad \text{for} \quad |z|<1.
\end{equation}
Fix any two distinct points $z,w \in \D$ and consider the function 
\[ u(\lambda)=\log \frac{|f^\lambda(z)-f^\lambda(w)|}{4}.
\]
In view of holomorphic dependence and injectivity of the flow $\{f^\lambda\}$ \cite[Theorem 5.7.3]{aim09} $u$ is a harmonic function. Moreover, $u$ is negative by \eqref{eq:koebe}. An application of Harnack's inequality then gives
\[ u\left( \frac{K-1}{K+1} \right) \leq \frac{1}{K} u(0).
\]
Observe that for $\lambda=(K-1)/(K+1)$ the solution of \eqref{eq:beltrami} is the original mapping $f$ and for $\lambda=0$ the solution is the identity map. Thus after exponentiating the previous inequality yields
\[ \frac{|f(z)-f(w)|}{4} \leq \left( \frac{|z-w|}{4} \right)^{1/K},
\]
as stated in the Theorem.

It remains to show sharpness of the constant $4^{1-1/K}$. The examples will be based on quasiconformal deformation of ellipses. Let $R>1$ and denote by $\E_R$ the ellipse given by
\[ |z-2|+|z+2| < 2\, (R+1/R).
\]
For any $K \ge 1$ there is a $K$-quasiconformal mapping $g \colon \E_R \to \E_{R^{1/K}}$ (onto) which fixes the foci $g(\pm2)=\pm2$ and have affine boundary values.

Indeed, it is easy to construct such a mapping explicitly.
The conformal mapping $\phi(z)=z+1/z$ maps the annulus $A_R=\{z \colon 1<|z|<R \}$ onto $\E_R \setminus [-2,2]$ and extends as an affine map to $\{z \colon |z| \ge R\}$. Now the desired map is given by $g=\phi \circ \rho \circ \phi^{-1}$, where $\rho(z)= z|z|^{1/K-1}$ is a $K$-quasiconformal radial stretching. Note that $g$ extends over the slit $[-2,2]$ and thus we have the required $K$-quasiconformal map $g \colon \E_R \to \E_{R^{1/K}}$.

Let  $\alpha_R$ be an affine map from $\D$ onto $\E_R$. Consider now the map $f:=\alpha_{R^{1/K}}^{-1} \circ g \circ \alpha_R$. This is a $K \cdot(1+\varepsilon)^2$-quasiconformal selfmap of the unit disk with identity boundary values (possibly, after a rotation) where $\varepsilon=\varepsilon(R^{1/K}) \to 0$ as $R \to \infty$. Furthermore, $f$ maps a segment of length $\sim \frac{4}{R}$ to a segment of length $\sim \frac{4}{R^{1/K}}$ as $R \to \infty$. Indeed, the segment $\alpha_R^{-1}([-2,2])$ is mapped onto $\alpha_{R^{1/K}}^{-1}([-2,2])$. This shows that the constant $4^{1-1/K}$ is best possible in \eqref{eq:main}.
\end{proof}

\bibliographystyle{amsplain}

\end{document}